\pdfoutput=1
\documentclass{article}
\usepackage[margin=25mm]{geometry}
\usepackage{amsmath}
\usepackage{amsfonts}
\usepackage{amssymb}
\usepackage{graphicx}

\usepackage[utf8]{inputenc} 
\usepackage[T2A]{fontenc}
\usepackage[english]{babel}
\usepackage[unicode, pdftex]{hyperref}

\usepackage{xcolor}
\usepackage[all]{xy}

\title{Regular octagons in Poincaré model of Lobachevsky geometry}
\author{Oleksandr Pryshliak}

\date{\today}

\begin{document}
\maketitle

\begin{abstract}
To investigate the topological structure of Morse flows on the 2-disk we use the planar graphs as destinguished graph of the flow. We assume, that the flow is transversal to the boundary of the 2-disk.
We give a list of all planar graph with at least 3 edges and describe all planar graphs with 4 edges. We use a list of spherical graph with at least 4 edges.
\end{abstract}
\textit{Key words and phrases.} Morse flow, planar graphs, spherical graph, topological invariant.

\subsection*{Introduction}

When building hyperbolic structures on closed surfaces, you can use Lobachevsky's geometry on the plane. For this, it is necessary to present the surface in the form of a 2n-gon on the Lobachevsky plane and set the action of a discrete group, which is a subgroup of motions of the Lobachevsky plane, for which the 2n-gon is the fundamental domain. If such a surface is a double torus (an oriented surface of genus 2), then it can be obtained by gluing opposite sides of an octagon. In fact, the Lobachevsky plane is divided into octagons. The presence of symmetries simplifies calculations, then the problem of dividing into regular octagons naturally arises. In addition, it is important to give examples of such octagons, setting the coordinates of their vertices in one of Lobachevsky's models of geometry.

\subsection{Relevance of the topic}

Lobachevsky's geometry is a relevant topic, as it is an important component of mathematics and has applications in various fields of science and technology.

First of all, Lobachevsky's geometry studies geometry on non-Euclidean spaces, which differs from traditional Euclidean geometry, which uses Euclid's axiom of parallelism. This allows solving problems that cannot be solved in ordinary geometry. Lobachevsky's geometry has important applications in mathematical physics, cosmology, number theory, group theory, and game theory.

In addition, Lobachevsky's geometry finds application in geodesy and cartography, where non-Euclidean geometries are used to develop global navigation systems and maps.

In addition, Lobachevsky's geometry is important in the field of education, where the study of non-Euclidean geometry contributes to the development of abstract thinking and broadening the horizons of students.

Lobachevsky's geometry is important not only from a mathematical point of view, but also from a practical point of view.

One of the important fields of application of Lobachevsky's geometry is computer graphics and image processing. Lobachevsky's geometry is used in the creation of three-dimensional objects, as well as in the modeling of various processes and systems. In addition, non-Euclidean geometries are important in geometric image processing and pattern recognition.

Lobachevsky's geometry is also important in topology, which studies the properties of geometric shapes that remain unchanged under continuous mapping. Non-Euclidean geometries are used in the study of topological spaces and forms.

It is also important to note that Lobachevsky's geometry is important in philosophy because it reflects different approaches to understanding space and geometric form. It allows you to study and compare different concepts of space, which helps to solve philosophical questions about human perception of the world.

\subsection{Literature review}

There are many books on Lobachevsky's geometry, both domestic and foreign

A. V. Pogorelov's book \cite{pogorelov2003} is a well-known and widespread book for studying non-Euclidean geometry. The book is intended for students of mathematical and physical specialties, as well as for teachers who teach this course.

In the book, the author offers readers an introduction to non-Euclidean geometry, in particular the geometry of Lobachevsky and Riemann. The author describes in detail the basic concepts and properties of non-Euclidean geometry, such as geometric shapes, the relationship between angles and sides, the distance between points, and others.

The book also contains a large number of examples and exercises that help readers understand and apply theoretical knowledge in practice. The author uses easy-to-understand language and well-known examples that help readers learn complex material.

Yu. B. Shevchenko's book \cite{shevchenko2001} contains introductory materials on Lobachevsky's geometry, which allows one to imagine non-Euclidean spaces and their properties. The book contains various examples, problems and exercises that help to understand the material and learn it.

The author also examines Lobachevsky's theory in contrast to classical Euclidean geometry, which allows the reader to understand the differences between them. The book is recommended for students of higher educational institutions studying mathematics and physics.

E. H. Gladkyi's book \cite{gladkyi2007} contains introductory materials on Lobachevsky's geometry, which allows the reader to understand non-Euclidean spaces and their properties. The book examines issues related to non-Euclidean spaces, in particular, characteristics of the geometry of such spaces.

The author offers various examples, problems and exercises that help to understand the material and learn it. The book is recommended for students of higher educational institutions studying mathematics and physics.

Yu. B. Shevchenko's book \cite{shevchenko2002} is an important source for studying non-Euclidean geometry, in particular Lobachevsky's geometry. The author describes in detail the main concepts and properties of this geometry, in particular, different types of geometric objects and their relationships in non-Euclidean spaces.

The book also contains a large number of examples, problems and exercises that help readers understand and apply theoretical knowledge in practice. This book is recommended for students and teachers of mathematics and physics at higher educational institutions.

The book by V. O. Goldstein \cite{goldstein2008} offers readers a systematic and logical approach to the study of Lobachevsky's geometry. The author offers various examples and illustrations that help to understand the material and learn it.

The book contains a detailed description of non-Euclidean spaces and their properties, in particular, the properties of geometric objects that differ from classical Euclidean geometry. This book is also recommended for students and teachers of mathematics and physics at higher educational institutions.

Roberto Bonola's book \cite{bonola} provides a historical overview of the development of non-Euclidean geometry, including the works of Gauss, Lobachevsky, Bolyai, Riemann, and others. The book also contains a detailed discussion of the axiomatic approach to geometry and the controversy surrounding the fifth postulate.

\cite{ryan} Patrick Ryan's book offers an introduction to Euclidean and non-Euclidean geometry using an analytical approach. The book contains a discussion of the axiomatic method, but the emphasis is on the use of analytical methods for the study of geometry. The book includes exercises and solutions.

Book \cite{sommerville} D.M.I. Summerville provides a comprehensive introduction to non-Euclidean geometry, including hyperbolic and elliptic geometry. The book includes a discussion of the axiomatic approach and the various models used to represent non-Euclidean geometries. The book also includes a chapter on the history of non-Euclidean geometry and its relationship to other branches of mathematics.

In the book \cite{hadamard} Jacques Hadamard, an outstanding mathematician of the 20th century, examines the basic concepts and theorems of hyperbolic geometry, in particular, the geometric structure of Minkowski space is studied, the importance and significance of hyperbolic geometry in various areas of mathematics and physics are emphasized.

\cite{curtiss}'s book is an introduction to non-Euclidean geometry, in particular Lobachevsky's geometry and Riemannian geometry. The author offers an analytical approach to non-Euclidean geometry that allows mathematicians and physicists to apply these ideas in practical problems. The book contains a large number of examples and exercises that help to understand the material.

The book \cite{canary} contains lectures by famous mathematicians devoted to hyperbolic geometry. The authors consider basic concepts and results, such as the geometric structure of hyperbolic space, geodesic lines, the concept of curves, the calculation of geometric quantities, etc. The book contains many illustrations and examples that help to understand complex theoretical concepts.

\cite{ratcliffe} is a high-quality book that covers the fundamental theory of hyperbolic manifolds in detail. The author shows how these manifolds are related to geometric structures of different dimensions, and investigates their properties. The book contains many examples and illustrations to help the reader understand complex concepts.

\cite{anderson}'s book is an introduction to hyperbolic geometry. The author describes the basic concepts and results of hyperbolic geometry, such as the Poincaré model and hyperbolic manifolds. The book includes many examples and exercises to help the reader understand the material.

\cite{wolfe}'s book is an introduction to Lobachevsky geometry and hyperbolic geometry in general. The author provides high-quality and accessible explanations of basic concepts and theorems related to hyperbolic geometry. The book contains many examples and illustrations that help to understand the material.

\cite{tabachnikov}'s book is an introduction to hyperbolic geometry, with a focus on group theory. The author investigates the geometric constructions that arise when studying the symmetry of hyperbolic spaces, and describes in detail the properties of these constructions. The book contains many illustrations and examples

The basic concepts of hyperbolic and Riemannian geometry and topology are also contained in  \cite{prishlyak2012topological, prish2002theory, prish2004difgeom, prish2006osnovy, prish2015top}.

If we specify the gluing of the sides of the polygon using the action of a discrete subgroup of the group of transformations of the Lobachevsky plane, then we obtain a surface and a graph embedded in it. Different gluings will lead to different graphs. Such objects naturally arise as topological invariants of functions and flows. In particular, topological invariants of functions were constructed in papers \cite{Kronrod1950} and \cite{Reeb1946} and in \cite{lychak2009morse} for unoriented surfaces, and in \cite{Bolsinov2004, hladysh2017topology, hladysh2019simple, prishlyak2012topological} for of surfaces with a boundary, in \cite{prishlyak2002morse} for non-compact surfaces. Graphs as topological invariants of smooth functions were also studied in papers \cite{bilun2023morseRP2, bilun2023morse, hladysh2019simple, hladysh2017topology, prishlyak2002morse, prishlyak2000conjugacy, prishlyak2007classification, lychak2009morse, prishlyak2002ms, prish2015top, prish1998sopr, bilun2002closed, Sharko1993}, for manifolds with a boundary in papers \cite {hladysh2016functions, hladysh2019simple, hladysh2020deformations}, and on 3- and 4-dimensional manifolds in \cite{prishlyak1999equivalence, prishlyak2001conjugacy}. Topological flow invariants were studied in \cite{bilun2023gradient, Kybalko2018, Oshemkov1998, Peixoto1973, prishlyak1997graphs, prishlyak2020three, akchurin2022three, prishlyak2022topological, prishlyak2017morse, kkp2013, prishlyak2021flows, prishlyak2020topology, prishlyak2019optimal, prishlyak2022Boy},\cite{bilun2023discrete, bilun2023typical, loseva2016topology, loseva2022topological, prishlyak2017morse, prishlyak2022topological, prishlyak2003sum, prishlyak2003topological, prishlyak1997graphs, prishlyak2019optimal, stas2023structures} and  \cite{prish1998vek, prish2001top, Prishlyak2002beh2, prishlyak2002ms, prishlyak2007complete, hatamian2020heegaard, bilun2022morse, bilun2022visualization}.

\subsection{The purpose of the paper}

The goal of the paper is to construct regular octagons (according to the coordinates of their vertices and the midpoints of the sides in the Poincaré model), which can be used to tile the entire Lobachevsky plane.



\section{History of emergence and basic concepts of Lobachevsky's geometry}
\subsection{Emergence History}

In the 3rd century BC, the Greek scientist Euclid systematized the geometric information known to him in the great work "The Beginning". For more than two thousand years, this book has served as a geometry textbook around the world.

Careful study of Euclid's system has led scientists to the conclusion that "Beginnings" has quite serious flaws. For example, the number of axioms formulated by Euclid is not enough for a rigorous presentation of geometry, so Euclid relied on direct evidence, visualization, intuition, and sensory perceptions when presenting some of his proofs.
In addition to the geometry studied at school (Euclid's geometry), there is another geometry, Lobachevsky's geometry. This geometry is significantly different from Euclid's, for example, it states that many lines can be drawn through this point, parallel to this line, that the sum of the angles of a triangle is less than 180o. In Lobachevsky's geometry, there are no rectangles, similar triangles, etc.

Non-Euclidean geometry arose as a result of long attempts to prove Euclid's V postulate - the axiom of parallelism. This geometry is in many ways surprising, unusual and in many ways does not correspond to our usual ideas about the real world. But logically, this geometry is not inferior to the geometry of Euclid.

In the history of the development of the axiomatic method, an important role was played by the axioms of D. Hilbert, a German scientist (1862-1943), who stood out among the constellation of scientists of that period. These axioms at one time corresponded to the level of rigor of geometry. In 1899, D. Hilbert wrote: "Geometry, like arithmetic, requires only a few simple basic propositions to be constructed. These basic theses are called axioms of geometry. Establishing the axioms of geometry and studying their relationships is a task that has been the subject of numerous beautiful works of mathematical literature since the time of Euclid. This task boils down to a logical analysis of our spatial representation.

The axiomatic method, first developed by D. Hilbert in geometry from new positions, penetrated into other branches of mathematics: set theory, algebra, topology, probability theory, etc. In addition, the axiomatic method began to be used in the construction of other sciences, especially physics. These achievements are connected with the revolution in geometry made by M.I. Lobachevsky. Historically, Euclid's fifth postulate has attracted the attention of mathematicians for many centuries. Having deeply analyzed attempts to prove the fifth postulate, both his own and those belonging to other mathematicians, M.I. Lobachevsky came to the conviction of the independence of this postulate from other axioms, that is, the non-contradiction of geometry, in which the existence of two different straight lines that pass through this point parallel to a given straight line is axiomized.
E. Lobachevsky not only predicted the existence of a new geometry - non-Euclidean, but also developed it in detail. His opinion contradicted all the ideas of man about the surrounding world. The new geometry sharply diverged from the philosophical view of space (I. Kant), so this discovery was stunning. It turned out that the assumption about the non-Euclidean nature of real physical space did not contradict Euclid's axioms, except for the fifth postulate.

In the 70s of the last century, the non-contradiction of geometry was proved, which rightfully received the name of Lobachevsky. This proof was constructed using the Cayley-Klein and Poincaré models.

Among Euclid's axioms was the axiom about the parallelism of straight lines, or more precisely, the fifth postulate about parallel lines: if two straight lines create internal angles with a third on one side of it, the sum of which is less than the extended angle, then such straight lines intersect if they are extended sufficiently on one side .

In the modern formulation, it means the existence of no more than one straight line passing through a given point outside a given straight line and parallel to this straight line.

The complexity of the formulation of the fifth postulate gave rise to the idea of its possible dependence on other postulates, and therefore there were attempts to derive it from other premises of geometry. All attempts ended in failure. There were attempts to prove from the opposite: to achieve wiping, assuming the prevention postulate to be true. However, this path was unsuccessful.
It turned out that the fifth postulate does not depend on the previous ones, and therefore, it can be replaced by an equivalent one. And at the beginning of the 19th century, almost simultaneously, several mathematicians: K. Gauss in Germany, J. Bolya in Hungary, and N. Lobachevsky in Russia, had the idea of the existence of a geometry in which the axiom that replaces the fifth is true postulate: on a plane, through a point that does not lie on a given line, there are at least two straight lines that do not cross the given line.

Due to the priority of M. Lobachevsky, who was the first to come up with this idea in 1826, and his contribution to the development of a new, different from Euclidean geometry, the latter was named "Lobachevsky's geometry" in his honor.
The axiom of Lobachevsky's planimetry differs from the axiom of Euclid's planimetry by only one axiom: the axiom of parallelism is replaced by its negation - Lobachevsky's axiom of parallelism: through a point that does not belong to a given line, in the plane defined by them, at least two non-intersecting lines can be drawn.

\subsection{Basic Concepts}

The plane (or space) in which the Lobachevsky axiom holds is called the Lobachevsky plane (space).
Let us consider the main consequences of Lobachevsky's axiom.

\textbf{Lemma}. If at the intersection of two straight lines the intersecting angles (corresponding angles) are equal, then the straight lines do not intersect.

The second statement of the theorem directly merges with the proven.
The question arises: how many straight lines parallel to the straight line a pass through the point M, which does not lie on the straight line a? The answer to this question is given by the following theorem.

\textbf{Theorem}. If postulate V holds, then only one line parallel to line a passes through each point M that does not lie on line a.

Let's draw an arbitrary straight line b' different from straight line b through point M. One of the adjacent angles 1 or 2, indicated in the same figure, is acute; let1 be sharp. At the intersection of straight lines a and b' with straight line MN, we obtain internal one-sided angles: 1 and 3, the sum of which is less than two right angles, therefore, according to the V postulate, straight lines a and b' intersect
There is also an inverse theorem:

\textbf{Theorem}. If we assume that only one parallel line passes through a point that does not lie on this line, then postulate V is valid.

Therefore, the V postulate is equivalent (equivalence) to the so-called axiom of parallel lines: no more than one parallel line passes through a point that does not lie on this line.

\textbf{Theorem}. There are three types of mutual placement of two straight lines on the Lobachevsky plane: coincident, divergent, parallel.

Let M be a point that does not lie on line a, and let MN be the perpendicular drawn from point M to line a. Let's choose two points A and B on the line so that A — N — B. It follows from the theorem that through point M there is a single straight line CD parallel to the straight line AB and a single straight line EF parallel to the directed straight line BA.
The angles DMN and FMN are acute, therefore CD and EF are different straight lines. Let's prove that DMN =FMN. Let, on the contrary, DMN=FMN. eg DMN > FMN. Consider the ray MF', symmetric MF with respect to the straight line MN. This ray is the inner ray of angle DMN. Since MF does not intersect the line AB, then MF' does not intersect this line. But this contradicts the definition of the parallelism of lines CD and AB.
Thus, through each point M that does not lie on this line a, two lines parallel to line a pass in two different directions. These straight lines form equal sharp lines with the perpendicular MN drawn from the point to the straight line a. Each of these angles is called an angle of parallelism at point M with respect to line a.

Thus, two lines parallel to line a pass through the point: one parallel in each of the two directions on line a.

Lines a and b are said to diverge if they lie in the same plane, do not intersect and are not parallel. There are many outstanding lines.

The value of the angle of parallelism is completely determined by the distance from point M to straight line a.

In Lobachevsky's geometry, there is a dependence between angular and linear quantities; this is the essential difference between Lobachevsky's geometry and Euclid's geometry. Thus, in Lobachevsky's geometry there is no similarity of figures; in particular, triangles with correspondingly equal angles are congruent. Another feature of Lobachevsky's geometry is related to the unit of measurement of lengths. In Euclidean geometry, there are absolute constants of angular quantities, such as the right angle or radian, while there are no absolute linear constants. In order to express the lengths of segments in numbers, it is necessary to choose a unit of length measurement. A random segment can be chosen as such a unit. In contrast, in Lobachevsky's geometry, there is no need, since, having a natural unit of measurement of angles, it is possible to agree on the choice of a natural unit of lengths. Example,

Let us now consider the simplest curves in Lobachevsky's geometry.

With the help of three types of beams, analogues of the circles of the Euclid plane are obtained. To do this, it is enough to consider the orthogonal trajectories of beams, that is, lines that cross all straight beams at right angles.

1. A circle is an orthogonal trajectory of straight lines of the first kind. Any size of the radius is possible.

2. Equidistant, or a line of equal distances - an orthogonal trajectory of straight lines of the second kind, the base is not included here. It is proved that these equidistant points are at a constant distance from the base. This line is concave towards the base and is not closed. In Klein's interpretation, the geometry of the Lobachevsky equidistant is considered as a circle with an ideal center, that is, with a center lying beyond infinitely distant points.

3. The boundary line, or oricycle, is an orthogonal trajectory of straight bundles of the third kind. It has excellent properties. First of all, oricycles are congruent. In general, two curves are called congruent, if between the points of these curves such a mutually unambiguous correspondence can be established, in which each chord of one curve is congruent to the corresponding (that is, connecting the corresponding points) chord of the other.

In addition, the oricycles are not closed and are concave in the direction of parallelism. The boundary line can be a circle with an infinitely distant center. At the same time, two boundary lines constructed for the same bundle cut congruent segments on the straight lines of this bundle: AA((BB(. Thus, the oricycles are concentric. And another property of them is that the ratio of the lengths of the arcs of the oricycles)), contained between two straight beams is an indicative function of the distance between the arcs.

\section{Poincaré model of Lobachevsky geometry}

The Poincaré model of Lobachevsky geometry was developed by the French mathematician Henri Poincaré in the 1880s. He was interested in the study of nonlinear differential equations and topology problems that arose in the context of Lobachevsky's geometry.

In his research, Poincaré applied the ideas of projective geometry, which allowed studying the properties of geometric objects while preserving their relative location. He proposed a model of Lobachevsky geometry that used the properties of projective geometry to represent a hyperbolic plane in three-dimensional space.

According to the Poincaré model, the hyperbolic plane can be represented as a circle inside a sphere called the Poincaré sphere. This sphere has a number of specific properties that allow representing geometric objects in the hyperbolic plane and performing various geometric operations.

The Poincaré model of Lobachevsky geometry is one of the most famous and widespread models of Lobachevsky geometry, which is used in various fields of mathematics, physics, and computer science.

Poincaré model of Lobachevsky geometry on the unit circle
$x^2+y^2<1$ is given by the Riemannian metric $$ds^2=4\frac{dx^2+dy^2}{(1-x^2-y^2)^2}.$$ 

The limit circle $x^2+y^2 =1$ is called an absolute for
Poincaré models.

Hyperbolic lines are called intersections of lines or circles
perpendicular to the absolute with the circle
Poincaré, respectively. The distance between points is found as length
segment of a hyperbolic line with endpoints at these points in the corresponding
Riemannian metric.

If the points are set using complex
numbers, then the distance between them can
be calculated by formulas

$$d(z,w)=\ln \frac{|1-z\overline{w}|+|z-w|}{|1-z\overline{w}|-|z-w|},$$

$$\tanh \frac{d(z,w)}{2}=\frac{z-w}{1-z\overline{w}}.$$

\section{Trigonometric relations in Lobachevsky's geometry}
\subsection{Properties of triangles in the Lobachevsky plane}

As Legendre showed, the statement that the sum of the angles of a triangle is equal to the postulate is not equivalent. Therefore, on the Lobachevsky plane, the sum of the angles of a triangle should not be equal. It cannot be greater than, which contradicts absolute geometry. Therefore, the sum of the angles of the triangle is smaller on the Lobachevsky plane.

The question arises: what number is the sum of the angles of a triangle in the Lobachevsky plane? It turns out that such a specific number for triangles in Lobachevsky's geometry does not exist. The following theorem is valid here.

\textbf{Theorem.} The sum of the angles of a triangle in Lobachevsky's geometry is not a constant value and depends on the shape and dimensions of the triangle.

In Euclid's geometry, the angular defect of a triangle is zero, and in Lobachevsky's geometry, the angular defect is a variable value that varies from zero to
At the same time, in Lobachevsky's geometry, the area of the triangle is proportional to the angular defect.

\textbf{Theorem}. The sum of the angles of any triangle is not greater than 2d.

Therefore, the sum of the angles of any triangle is not greater than 2d. But can it not turn out that in some triangles this sum is less than 2d, and in others it is equal to 2d? The negative answer to this question is given by the second theorem of Saccera-Legendre.

\textbf{Theorem}. If the sum of the angles of one triangle is equal to 2d, then the sum of the angles of any triangle is equal to 2d.
So, we get another assumption equivalent to postulate V: there is at least one triangle whose sum of angles is equal to 2d.

\textbf{Theorem.}
If on the Lobachevsky plane three angles of one triangle are equal to three angles of another triangle, then such triangles are congruent.

This theorem gives a new, fourth sign of the equality of triangles, which is not found in absolute geometry.

\textbf{Theorem.} On the Lobachevsky plane, a circle can be described around any triangle.

\textbf{Theorem.} If two median perpendiculars to the sides of a triangle intersect, then the third median perpendicular also passes through this point of intersection.

\textbf{Corollary.} If two median perpendiculars to the sides of a triangle intersect at point O, then a circle centered at point O can be described around such a triangle.

\textbf{Theorem}. If the two median perpendiculars to the sides of the triangle diverge, then the median perpendicular to the third side also diverges with the first two and all have a single common perpendicular, and all the vertices of the triangle are equidistant from it.

\textbf{Corollary.} If the average perpendiculars to two sides of a triangle diverge, an equidistant can be described around the triangle.

\textbf{Theorem.} If two median perpendiculars oriented in the same direction to the sides of the triangle are parallel, then the median perpendicular to the third side of the triangle is parallel to the first two.

Proof. The validity of this theorem follows from the previous two. Indeed, if we assume that the median perpendicular from the AC side intersects with the median perpendicular a to the AB side, then according to Theorem 6, the third perpendicular b must also pass through the point of their intersection, and according to Theorem 7, it must diverge, which contradicts the condition of the theorem. So, the middle perpendicular with parallel a and b.

\textbf{Corollary.} If the median perpendiculars to two sides of a triangle are parallel in this direction, an oricycle can be described around such a triangle.
It follows from these theorems that one of three lines can be described in the Lobachevsky plane of each triangle - either a circle, an equidistant, or an oricycle.

\subsection{Theorems of sines and cosines}

In Lobachevsky geometry, for any triangle with sides $a$, $b$ and $c$ and corresponding opposite angles $\alpha$, $\beta$ and $\gamma$, the theorem of sines has the following form:

$$\frac{\text{sinh} a}{\text{sin} \alpha} =
\frac{\text{sinh} b}{\text{sin} \beta} =
\frac{\text{sinh} c}{\text{sin} \gamma} $$

where $\sinh$ is the hyperbolic sine and $\sin$ is the trigonometric sine.

In Lobachevsky geometry, for any triangle with sides $a$, $b$ and $c$ and corresponding opposite angles $\alpha$, $\beta$ and $\gamma$, the first theorem of cosines has the following form:

$$\cosh c =\cosh a \cosh b + \sinh a \sinh b \cos \gamma$$

where $\cosh$ is the hyperbolic cosine, $\sinh$ is the hyperbolic sine, and $\cos$ is the trigonometric cosine.

In Lobachevsky geometry, for any triangle with sides $a$, $b$ and $c$ and corresponding opposite angles $\alpha$, $\beta$ and $\gamma$, the second theorem of cosines has the following form:

$$\cos \alpha = -\cos \beta \cos \gamma + \sin \beta \sin \gamma \cosh a$$

where $\cosh$ is the hyperbolic cosine, $\sin$ is the trigonometric sine, and $\cos$ is the trigonometric cosine.

The last formula allows you to find the sides of a triangle by their angles.

\subsection{Right triangles}

For a right-angled triangle $ABC$ with a right angle $B$, the Pythagorean theorem holds:
$$\cosh b = \cosh a \cosh c.$$

\begin{figure}[ht!]
\center{\includegraphics[width=0.45\linewidth]{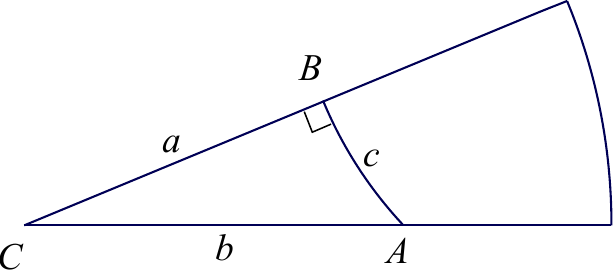}
}
\caption{Right triangle in the Poincaré model}
\label{triangle}
\end{figure}

If a regular octagon is cut along all axes of symmetry (straight lines passing through opposite corners and straight lines passing through the centers of opposite sides), we will get 16 right triangles, one of which is shown in fig. \ref{triangle}.

In this triangle, $\beta = \pi /2, \ \gamma = \pi /8$.

From the second theorem of cosines we have:

$$\cosh a = \frac{\cos \alpha }{\sin \gamma} = \frac{\cos \alpha }{\sin \pi /8}.$$

$$\cosh b = \frac{\cos \alpha \cos \gamma}{\sin \alpha \sin \gamma} = \frac{\cos \alpha \cos \pi /8}{\sin \alpha \sin \pi /8}.$$

\section{Regular octagons in the Poincaré model}

\begin{figure}[ht!]
\center{\includegraphics[width=0.60\linewidth]{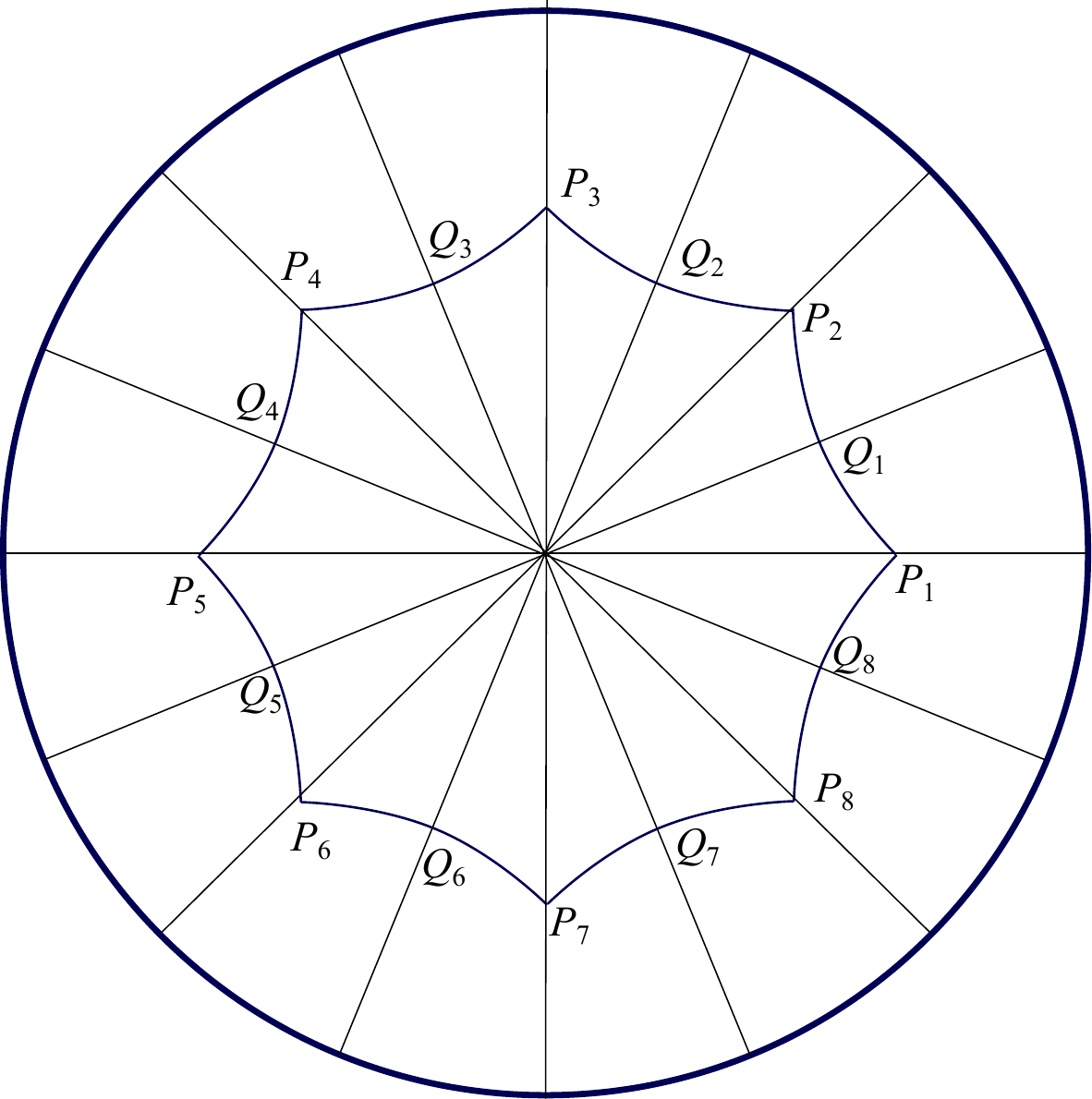}
}
\caption{Denotation of points on a regular octagon
}
\label{points}
\end{figure}

Denote the vertices of the regular octagon $P_1,P_2,P_3,P_4,$ $P_5,$ $P_6,$ $P_7,P_8$, and the centers of the sides $Q_1,Q_2,Q_3,$ $Q_4,$ $Q_5,$ $Q_6, Q_7,Q_8$ (see Fig. \ref{points}).

\subsection{Coordinates of vertices}

We will consider the Poincaré model in the unit circle $x^2+y^2<1$. $O (0,0)$ is the origin of coordinates.

Let the vertices have coordinates $P_i(x_i,y_i)$, and the midpoints of the sides $Q_i(u_i,v_i)$.

Let $\alpha$ be equal to half of the internal angle of the octagon.

Let's find the coordinates of the first vertex $P_1$. The triangle $P_1Q_1O$ coincides with the right triangle $ABC$ considered in the previous section. Therefore $$\cosh b = \cosh |OP_1| = \frac{\cos \alpha \cos \pi /8}{\sin \alpha \sin \pi /8}=\cot \alpha \cot \pi /8 .$$

On the other hand, $OP_1$ is equal to the length of the line segment $x=t, y=0, t\in [0,x_1]$. According to the arc length formula in the Poincaré model, we have
$$b=|OP_1|=\int_0^{x_1} \frac{2 dx}{1-x^2}=\log \frac{1+x_1}{1-x_1}.$$
We will get from here
$$x_1=\frac{e^b-1}{e^b+1}.$$

Substituting one formula into another, we have

$$x_1=\frac{e^{\cosh ^{-1}\left(\cot \left(\frac{\pi }{8}\right) \cot (\alpha)\right)}-1} {e^{\cosh ^{-1}\left(\cot \left(\frac{\pi }{8}\right) \cot (\alpha)\right)}+1}$$

Since $OP_2=OP_1$, and the angle $P_1OP_2$ is equal to $Pi/4$, then
$$x_2=y_2=\frac{x_1}{\sqrt{2}}.$$

Given the symmetry of a regular octagon, we find the coordinates of the remaining vertices:

$$P_3(0,x_1), \ P_4(-x_2,x_2),\ P_5(-x_1,0), \ P_6(-x_2,-x_2), \ P_7(0, -x_1), \ P_8(x_2 ,-x_2).$$

We will perform similar calculations for the coordinates of the midpoints of the sides (denote $u=u_1,v=v_1$):

$$\cosh a = |OQ_1| = \frac{\cos \alpha }{\sin \pi /8}.$$

Then $$ u=\frac{e^a-1}{e^a+1} \cos\frac{\pi}{8}, \
v=\frac{e^a-1}{e^a+1} \sin\frac{\pi}{8}.$$

The remaining coordinates of the midpoints of the sides are from symmetries:

$$Q_2(v,u), \ Q_3(-v,u),\ Q_4(-u,v),\ Q_5(-u,-v),$$ $$Q_6(-v,-u), \ Q_7(v,-u),\ Q_8(u,-v).$$

\subsection{Visualization of octagons}

To visualize octagons, let's draw their sides - these are circles passing through 3 points: $P_1Q_1P_2$, $P_2Q_2P_3$, etc.

To do this, we will use the formula for finding the center of a circle passing through three points $A(x_1, y_1)$, $B(x_2, y_2)$ and $C(x_3, y_3)$

$$x_0 = \frac{(y_2 - y_1)(y_3 - y_1)(y_3 - y_2) + x_1(y_2 - y_3)^2 + x_2(y_3 - y_1)^2 + x_3(y_1 - y_2)^2} {2(x_2 - x_1)(y_3 - y_1) - 2(x_3 - x_1)(y_2 - y_1)}$$

$$y_0 = \frac{(x_2 - x_1)(x_3 - x_1)(x_3 - x_2) + y_1(x_3 - x_2)^2 + y_2(x_1 - x_3)^2 + y_3(x_2 - x_1)^2} {2(y_2 - y_1)(x_3 - x_1) - 2(y_3 - y_1)(x_2 - x_1)}$$

The radius of a circle is equal to the distance from the center to any of the three points.

So, you need to know the coordinates of the first three points (the rest of the circles and points are found using symmetries relative to the coordinate axes or bisectors between them).

\textbf{Case} $\alpha = \pi/3$.

The results of calculations according to the above formulas:

$x_1= 0.405616$

$y_1=0$

$u= 0.336816$

$v= 0.139514$

$x_2= 0.286814$

$y_2= 0.286814$

\begin{figure}[ht!]
\center{\includegraphics[width=0.45\linewidth]{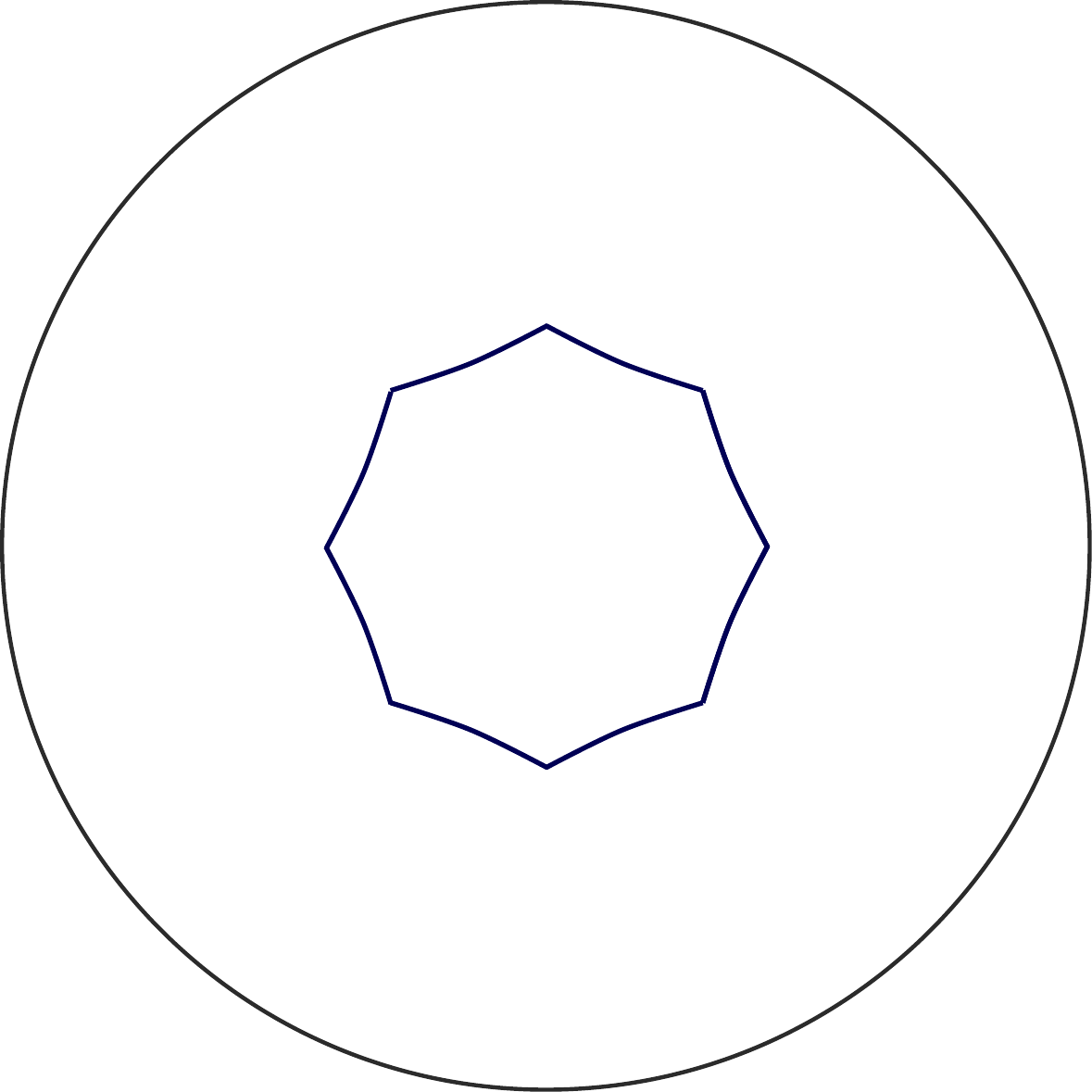}
}
\caption{A regular octagon with interior angle $2\alpha = 2\pi/3$}
\label{p3}
\end{figure}

\textbf{Case} $\alpha = \pi/4$.

The results of calculations according to the above formulas:

$x_1= 0.643594$

$y_1=0$

$u= 0.504081$

$v= 0.208797$

$x_2= 0.45509$

$y_2= 0.45509$

\begin{figure}[ht!]
\center{\includegraphics[width=0.45\linewidth]{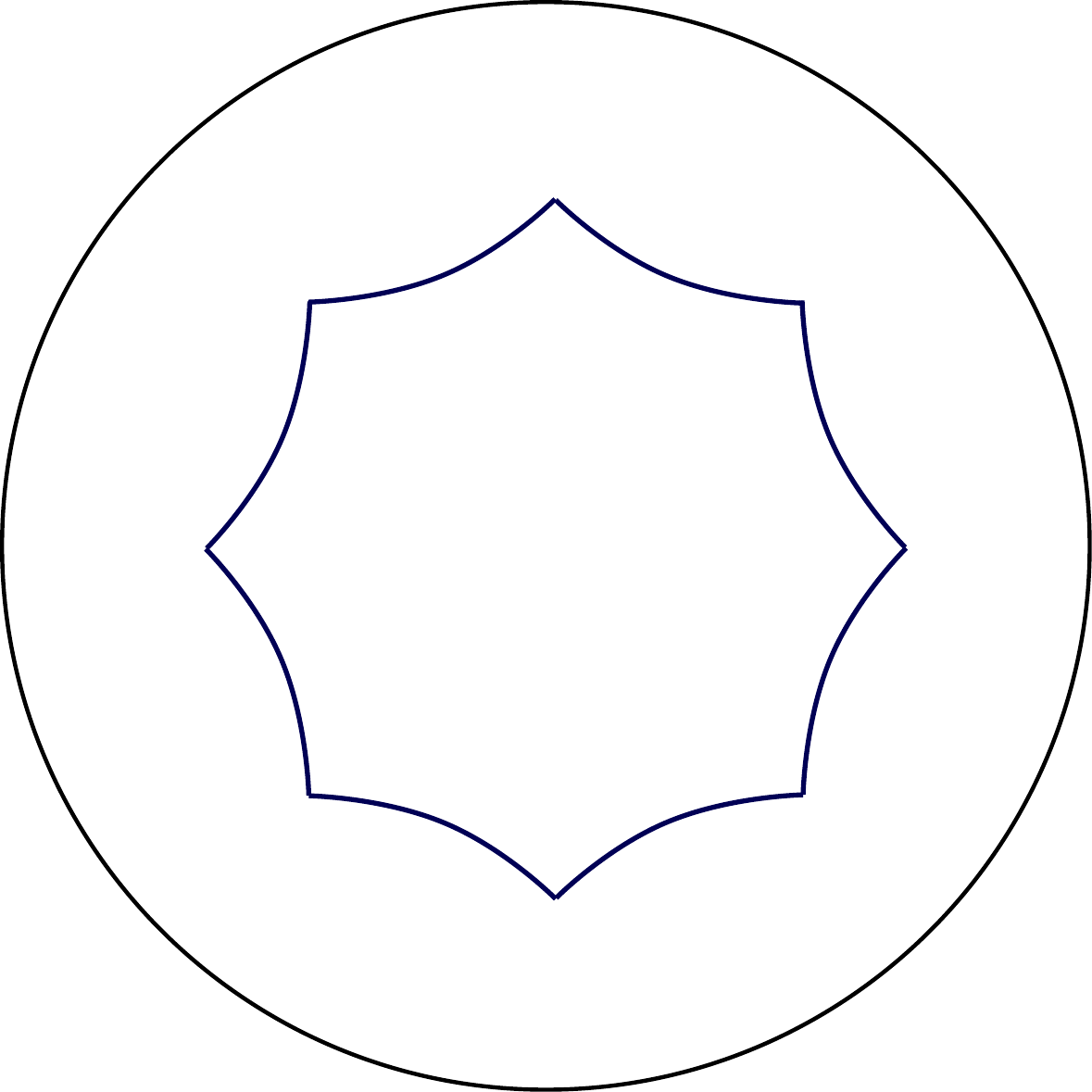}
}
\caption{Regular octagon with interior angle $2\alpha = \pi/2$}
\label{p4}
\end{figure}

\textbf{Case} $\alpha = \pi/6$.

The results of calculations according to the above formulas:

$x_1= 0.783591$

$y_1=0$

$u= 0.574794$

$v= 0.238087$

$x_2= 0.554082$

$y_2= 0.554082$

\begin{figure}[ht!]
\center{\includegraphics[width=0.45\linewidth]{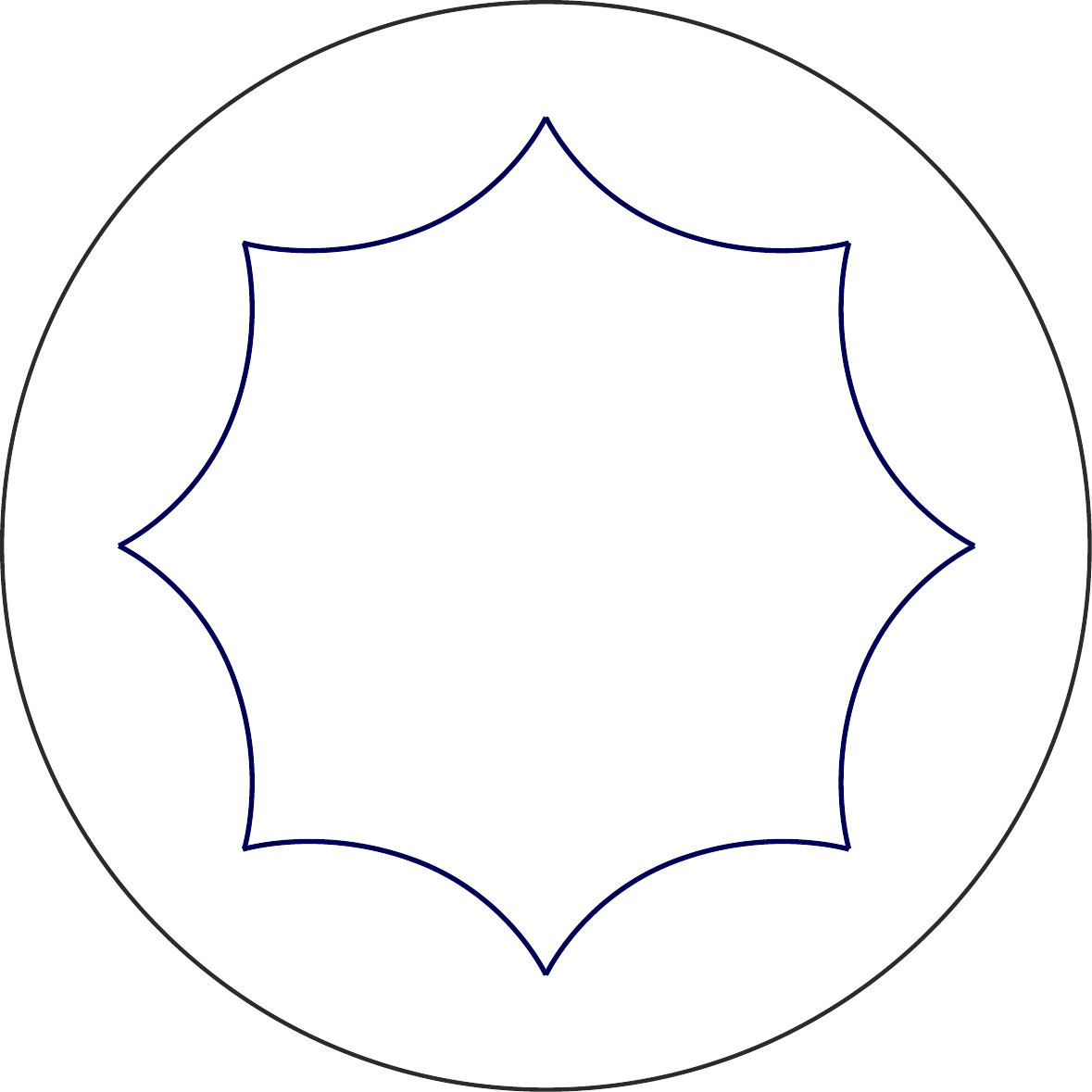}
}
\caption{Regular octagon with interior angle $2\alpha = \pi/3$}
\label{p6}
\end{figure}

\textbf{Case} $\alpha = \pi/8$.

The results of calculations according to the above formulas:

$x_1= 0.840896$

$y_1=0$

$u= 0.594604$

$v= 0.246293$

$x_2= 0.594604$

$y_2= 0.594604$

\begin{figure}[ht!]
\center{\includegraphics[width=0.45\linewidth]{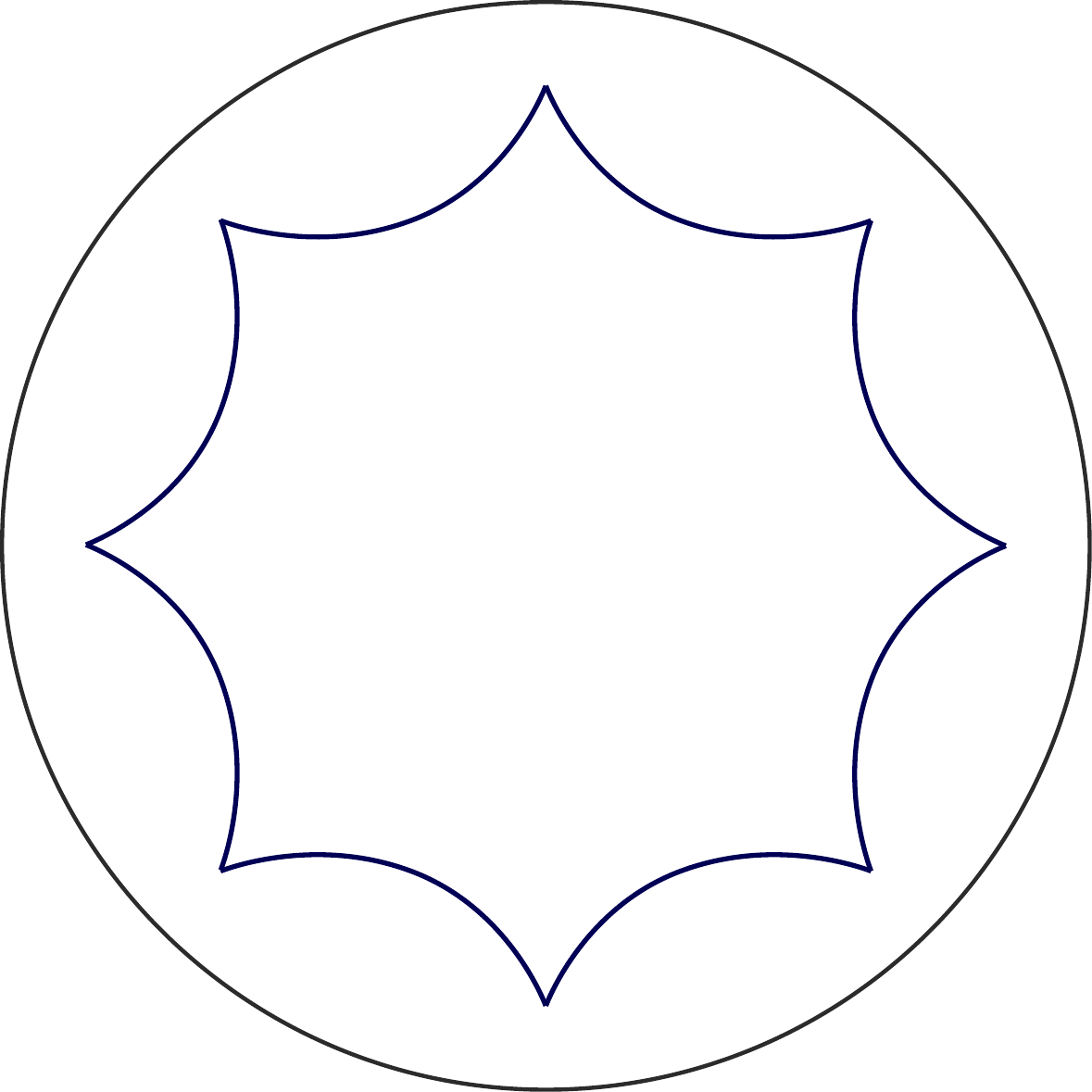}
}
\caption{Regular octagon with interior angle $2\alpha = \pi/4$}
\label{p8}
\end{figure}

\section*{Conclusion}
\addcontentsline{toc}{chapter}{Conclusion}

In the Poincaré model of Lobachevsky's geometry, regular octagons with different interior angles were constructed. The lengths of their diagonals and the coordinates of the vertices are found, provided that the center of the octagon is at the origin of the coordinates, and one of the vertices is on the abscissa axis. Such octagons were visualized for internal corners $2 \pi/3$, $\pi/2$,$\pi/3$,$\pi/4$.

The obtained results can be used to study hyperbolic and analytic structures on the double torus.

It would be interesting to generalize the obtained results to other surfaces and other regular polygons in Lobachevsky geometry.


\end{document}